
\magnification=1100
 \baselineskip=16truept

\footline={\tt\hss-\hskip-4pt-$
          \,$\folio$\,$-\hskip-4pt-\hss}

\font\title=cmbx12 scaled 1500
\font\section=cmbx12 scaled 1250
\font\subsection=cmcsc10 scaled\magstep1
\font\author=cmbx12 scaled 600

\font\tenmath=msbm10
\font\sevenmath=msbm7
\font\fivemath=msbm5
\newfam\mathfam \textfont\mathfam=\tenmath
\scriptfont\mathfam=\sevenmath
\scriptscriptfont\mathfam=\fivemath
\def\math{\fam\mathfam}
\def\r{{\math R}}
\def\e{{\math E}}
\def\p{{\math P}}


\font\msamten=msam10
\font\msamseven=msam7
\font\msamfive=msam5
\newfam\msamfam
\textfont\msamfam=\msamten
\scriptfont\msamfam=\msamseven
\scriptscriptfont\msamfam=\msamfive
\def\hexa #1{\ifcase #1 0\or 1\or 2\or 3\or 4\or 5\or 6\or 7\or
                        8\or 9\or A\or B\or C\or D\or E\or F\fi}

\def\qed{$\hfill\sqcup\!\!\!\!\sqcap$}
\def\ind{{\bf 1}\hskip-2.5pt{\rm l}}

\def\d{\, \hbox{\rm d}}



\vglue50pt

\centerline{\title On the increments of the}
\bigskip

\centerline{\title principal value of Brownian local time}
\bigskip
\bigskip
\bigskip
\bigskip
\bigskip

\noindent {\bf Endre C{\author S\'AKI}\footnote{$^1$}{\sevenrm Research
supported by the Hungarian National Foundation for Scientific Research,
Grant No. T 037886 and T 043037.}}
\vskip1ex

\noindent {\it Alfr\'ed R\'enyi Institute of Mathematics, Hungarian
Academy of Sciences, Budapest, P.O.B. 127, H-1364, Hungary.} {\tt E-mail:
\hskip-7pt csaki@renyi.hu}
\bigskip
\medskip

\noindent  {\bf Yueyun H{\author U}}
\vskip1ex

\noindent  {\it D\'epartement de Math\'ematiques, Institut
Galil\'ee (L.A.G.A. UMR 7539) Universit\'e Paris XIII,  99 Avenue
J-B Cl\'ement, 93430 Villetaneuse, France.} {\tt E-mail:
\hskip-7pt yueyun@math.univ-paris13.fr}

\bigskip
\bigskip

\noindent{\bf Summary.} Let $W$ be a one-dimensional Brownian
motion starting from 0. Define $Y(t)= \int_0^t{ \d s \over W(s)}
:= \lim_{\epsilon\to0} \int_0^t  1_{( |W(s)|> \epsilon)} { \d s
\over W(s)} $ as Cauchy's principal value related to local time.
We prove limsup and liminf results for the increments of $Y$.

\bigskip
\medskip
\noindent{\bf Running title.} Principal value increments.

\medskip


\medskip


\medskip

\vfill\eject

\baselineskip=18pt

\noindent  {\section 1. Introduction}
\bigskip

Let $\{ W(t); \, t\ge 0\}$ be a one-dimensional standard Brownian motion
with $W(0)=0$, and let $\{ L(t,x); \, t\ge 0, \, x\in \r\}$ denote its
jointly continuous local time process. That is, for any Borel function
$f\ge 0$,
$$
\int_0^t f(W(s))\d s= \int_{-\infty}^\infty f(x) L(t,x)\d x, \qquad t\ge 0.
$$

\noindent We are interested in the process
$$
Y(t):= \int_0^t {\! \d s\over W(s)}, \qquad t\ge 0. \leqno(1.1)
$$

\noindent Rigorously speaking, the integral $\int_0^t \d s/W(s)$ should
be considered in the sense of Cauchy's principal value, i.e., $Y(t)$ is
defined by
$$
Y(t):= \lim_{\varepsilon\to 0^+} \int_0^t {\! \d s\over W(s)}
\ind_{ \{ |W(s)|\ge \varepsilon\} }= \int_0^\infty \;
{L(t,x)-L(t,-x) \over x} \d x. \leqno(1.2)
$$

\noindent Since $x\mapsto L(t,x)$ is H\"older continuous of order $\nu$,
for any $\nu<1/2$, the integral on the extreme right in (1.2) is almost
surely absolutely convergent for all $t>0$. The process $\{Y(t),\, t\ge
0\}$ is called the principal value of Brownian local time.

It is easily seen that $Y(\cdot)$ inherits a scaling property from
Brownian motion, namely, for any fixed $a>0$, $t\mapsto
a^{-1/2}Y(at)$ has the same law as $t\mapsto Y(t)$. Although some
properties distinguish $Y(\cdot)$ from Brownian motion (in
particular, $Y(\cdot)$ is not a semimartingale), it is a kind of
folklore that $Y$ behaves somewhat like a Brownian motion. For
detailed studies and surveys on principal value, and relation to
Hilbert transform see Biane and Yor [4], Fitzsimmons and Getoor
[13],  Bertoin [2], [3],  Yamada [20], Boufoussi et al. [5], Ait
Ouahra and Eddahbi [1],  Cs\'aki et al. [11]
 and a collection of
papers [22] together with their references.  Biane and Yor [4]
presented a detailed study on $Y$ and determined a number of
distributions for principal values and related processes.

Concerning almost sure limit theorems for $Y$ and its increments,
we summarize the relevant results in the literature. It was shown
in  [17] that the following law of the iterated logarithm holds:

\medskip\noindent
{\bf Theorem A.} (Hu and Shi [17])
$$
 \limsup_{T\to \infty} \; {Y(T)\over \sqrt{T\log\log T}}
 = \sqrt{8}\, , \qquad \hbox{\rm a.s.} \leqno(1.3)
$$

This was extended in [10] to a Strassen-type [18] functional law of the
iterated logarithm.

\medskip\noindent
{\bf Theorem B.} (Cs\'aki et al. [10]) With probability one the set
$$
 \left\{{Y(xT)\over \sqrt{8T\log\log T}},\, 0\leq x\leq 1\right\}_{T\geq
3}\leqno(1.4)
$$
is relatively compact in $C[0,1]$ with limit set equal to
$$
{\cal S}:=
\left\{f\in C[0,1]:\, f(0)=0,\, f\,\, {\rm is\,\, absolutely\,\,
continuous\,\, and} \int_0^1(f'(x))^2\d x\leq 1\right\}.\leqno(1.5)
$$

\medskip
Concerning Chung-type law of the iterated logarithm, we have the
following result:

\medskip\noindent
{\bf Theorem C.} (Hu [16])
$$
\liminf_{T\to\infty}\sqrt{\log\log T\over T}\sup_{0\le s\le T}
|Y(s)|=K_1\, ,\qquad\hbox{\rm a.s.}\leqno(1.6)
$$
{\sl with some (unknown) constant $K_1>0$.}
\medskip

The large increments were studied in [7] and [8]:

 \noindent
{\bf Theorem D.} (Cs\'aki et al. [7]) {\sl Under the conditions
$$
\cases{
0<a_T\le T, \cr \cr
\hbox{$T\mapsto a_T$ {\rm and} $T\mapsto T/a_T$ {\rm are both
non-decreasing},} \cr\cr
{\displaystyle \lim_{T\to\infty}\; {\log (T/a_T)\over \log\log T}}=
\infty,
\cr}
\leqno(1.7)
$$
we have $$ \lim_{T\to\infty}\; {\sup_{0\le t\le T-a_T}\;
\sup_{0\le s\le a_T} \; |Y(t+s)-Y(t)|\over \sqrt{a_T \log
(T/a_T)}} =2, \qquad \hbox{\rm a.s.} \leqno(1.8)
$$ }

Wen [19] studied the lag increments of $Y$ and among others proved
the following results.
\noindent{\bf Theorem E.} (Wen [19]) {\sl }
$$
\limsup_{T\to\infty}\sup_{0\le t\le T}{\sup_{t\le s\le T}
|Y(s)-Y(s-t)|\over\sqrt{t(\log(T/t)+2\log\log t)}}=2,
\qquad\hbox{\rm a.s.}\leqno(1.9)
$$
{\sl Under the conditions $0<a_T\leq T$, $a_T\to\infty$ as
$T\to\infty$, we have}
$$
\limsup_{T\to\infty}\sup_{0\le t\le T-a_T}{\sup_{0\le s\le a_T}
|Y(t+s)-Y(t)|\over\sqrt{a_T(\log((t+a_T)/a_T)+2\log\log a_T)}}\leq 2,
\qquad\hbox{\rm a.s.}\leqno(1.10)
$$
{\sl If $a_T$ is onto, then we have equality in} (1.10).

\medskip
In this note our aim is to investigate further limsup and liminf behaviors
of the increments of $Y$.

\bigskip
\proclaim
Theorem 1.1.
Assume that $T\mapsto a_T$ is a function such that $0<a_T\le
T$, and both $a_T$ and $T/a_T$ are non-decreasing. Then
\par\noindent
{\rm (i)}
$$
\limsup_{T\to\infty}\; {\sup_{0\le t\le T-a_T}\; \sup_{0\le s\le a_T} \;
|Y(t+s)-Y(t)|\over \sqrt{a_T
\left(\log \sqrt {T/a_T}+\log\log T\right)}} =\sqrt{8}, \qquad
\hbox{\rm a.s.}\leqno(1.11)
$$
{\rm (iia)\quad} {\sl If $a_T>T(\log T)^{-\alpha}$ for some $\alpha
<2$, then}
$$
\liminf_{T\to\infty}\sqrt{\log\log T\over a_T}
\sup_{0\le t\le T-a_T}\sup_{0\le s\le a_T}|Y(t+s)-Y(t)|=K_2,
\qquad
\hbox{\rm a.s.}\leqno(1.12)
$$
{\rm (iib)\quad} {\sl If $a_T\le T(\log T)^{-\alpha}$ for some
$\alpha>2$,
then}
$$
\liminf_{T\to\infty}
{\sup_{0\le t\le T-a_T}\sup_{0\le s\le a_T}|Y(t+s)-Y(t)|
\over\sqrt{a_T\log(T/a_T)}}=K_3,
\qquad
\hbox{\rm a.s.}\leqno(1.13)
$$
{\sl with some positive constants $K_2, K_3$. If, moreover,}
$$
\lim_{T\to\infty}{\log(T/a_T)\over\log\log T}=\infty,
$$
{\sl then} $K_3=2$.

\proclaim
Theorem 1.2.
Assume that $T\mapsto a_T$ is a function such that $0<a_T\le
T$, and both $a_T$ and $T/a_T$ are non-decreasing. Then
\par\noindent
{\rm (i)}
$$
\liminf_{T\to\infty}\; { \sqrt{ T \log\log T} \over a_T}
\inf_{0\le t\le T-a_T}\; \sup_{0\le s\le a_T} \; |Y(t+s)-Y(t)| =K_4,
\qquad \hbox{\rm a.s.}\leqno(1.14)
$$
{\sl with some positive constant $K_4$.
If, $\lim_{T\to\infty}(a_T/T)=0$, then} $K_4=1/\sqrt{2}$.
\par\noindent
{\rm (iia)\quad} {\sl If} $0<\rho\leq 1$, {\sl then}
$$
\limsup_{T\to\infty}
{\inf_{0\le t\le T-\rho T}\; \sup_{0\le s\le \rho T} \; |Y(t+s)-Y(t)|
\over\sqrt{T\log\log T}} =\rho\sqrt{8},
\qquad \hbox{\rm a.s.}\leqno(1.15)
$$
{\rm (iib)\quad} {\sl If}
$$
\lim_{T\to\infty}{a_T(\log\log T)^2\over T}=0,
$$
{\sl then}
$$
\limsup_{T\to\infty}\; { \sqrt{T}  \over a_T \sqrt{\log\log T}}
\inf_{0\le t\le T-a_T}\; \sup_{0\le s\le a_T} \; |Y(t+s)-Y(t)| =K_5,
\qquad \hbox{\rm a.s.}\leqno(1.16)
$$
{\sl with some positive constant} $K_5$.

\bigskip\noindent
\proclaim Remark 1. The exact values of the constants $K_i$,
$i=2,3,4,5$ are unknown. It seems difficult to determine the exact
values of these constants. In the proofs we establish upper and
lower bounds with possibly different constants. It follows however
by 0-1 law for Brownian motion that the limsup's and liminf's
considered here are non-random constants.

\noindent \proclaim Remark 2.  Plainly we recover some previous
results on the path properties of $Y$  by considering particular
cases of Theorems 1.1 and 1.2. For instance, Theorems A and C
follow from (1.11) and (1.12) respectively by taking $a_T=T$, and
(1.8) follows from (1.11) combining with (1.13).  However in
Theorem 1.1(ii) and Theorem 1.2(ii) there are still  small gaps in
$a_T$.

The organization of the paper is as follows: In Section 2 some facts are
presented needed in the proofs. Section 3 contains the necessary
probability estimates. Theorem 1.1(i) and Theorem 1.1(iia,b) are
proved in Sections 4 and 5, resp., while Theorem 1.2(i) and Theorem
1.2(iia,b) are proved in Sections 6 and 7, resp.

  Throughout the paper, the letter $K$
with subscripts will denote some important but unknown finite
positive constants, while the letter $c$ with subscripts denotes
some finite and positive universal constants not important in our
investigations. When the constants depend on a parameter, say
$\delta$, they are denoted by $c(\delta)$ with subscripts.


\bigskip
\noindent
{\section 2. Facts}

Let $\{W(t),\, t\geq 0\}$ be a standard Brownian motion and define
the following objects: $$\leqalignno{
      g & :=\sup\{t:\, t\leq 1,\, W(t)=0\} &(2.1)
      \cr B(s)  &:={W(sg)\over \sqrt{g}},\qquad 0\leq s\leq 1, &(2.2)
\cr m(s)  &:={|W(g+s(1-g))|\over\sqrt{1-g}},\qquad 0\leq s\leq 1.
   &(2.3)
\cr}$$ Here we summarize some well-known facts needed in our
proofs.

\proclaim Fact 2.1. {\rm (Biane and Yor [4])}
$$
{\p(Y(1)\in \d x)\over\d x}=\sqrt{2\over
\pi^3}\sum_{k=0}^\infty(-1)^k \exp\left(-{(2k+1)^2x^2\over
8}\right), \quad x\in \r. \leqno(2.4)
$$
Consequently we have the estimate: for $\delta>0$
$$
c_1\exp\left(-{z^2\over 8(1-\delta)}\right)
\leq\p(Y(1)\geq z)\leq\exp\left(-{z^2\over 8}\right),
\qquad z\geq 1
\leqno(2.5)
$$
with some positive constant $c_1=c_1(\delta)$. Moreover, $g$,
$\{B(s),\, 0\leq s\leq 1\}$ and $\{m(s),\, 0\leq s\leq 1\}$ are
independent, $g$ has arcsine distribution, $B$ is a Brownian
bridge and $m$ is a Brownian meander.
$$
\eqalign{
& \p\left(\int_0^1{\d v\over m(v)}<z\, \Big|\, m(1)=0\right)
\cr & =\sum_{k=-\infty}^\infty (1-k^2z^2)\exp\left(-{k^2z^2\over 2}\right)
={8\pi^2\sqrt{2\pi}\over z^3}\sum_{k=1}^\infty
\exp\left(-{2k^2\pi^2\over z^2}\right),\quad z>0.\cr}
\leqno(2.6)
$$
$$
\p(m(1)>x)=  e^{-x^2/2},\qquad x>0.\leqno(2.7)
$$

\proclaim Fact 2.2. {\rm (Yor [21, Exercise 3.4 and pp. 44])}  Let
$Q_{x\to 0}^{\delta}$ be the law of square of Bessel bridge from
$x$ to $0$ of dimension $\delta>0$ during time interval $[0, 1]$.
The process $(m^2(1-v), 0\le v \le 1)$ conditioned on $\{
m^2(1)=x\}$ is distributed as $Q_{x\to 0}^3$. Furthermore, we have
$$
Q_{x\to 0}^{\delta}=Q_{0\to 0}^{\delta}*  Q_{x\to 0}^0,
 \quad \forall\, \delta>0, \, x >0, \leqno(2.8)
$$
where $*$ denotes convolution operator. Consequently, for any
$x>0$
$$
\p\left(\int_0^1{\d v\over m(v)}<z\, \Big|\, m(1)=x\right)\geq
\p\left(\int_0^1{\d v\over m(v)}<z\, \Big|\, m(1)=0\right).
\leqno(2.9)
$$

\proclaim Fact 2.3. {\rm (Hu [16])}
For $0<z\leq 1$
$$
c_2\exp\left(-{c_3\over z^2}\right)\leq \p(\sup_{0\leq s\leq
1}|Y(s)|<z)\leq c_4\exp\left(-{c_5\over z^2}\right)
\leqno(2.10)
$$
with some positive constants $c_2,c_3,c_4,c_5$.

\proclaim Fact 2.4. {\rm (Cs\"org\H o and R\'ev\'esz [12])}
Assume that $T\mapsto a_T$ is a function such that $0<a_T\le
T$, and both $a_T$ and $T/a_T$ are non-decreasing. Then
$$
\limsup_{T\to\infty}\; {\sup_{0\le t\le T-a_T}\; \sup_{0\le s\le a_T} \;
|W(t+s)-W(t)|\over \sqrt{a_T
\left(\log (T/a_T)+\log\log T\right)}} =\sqrt{2}, \qquad
\hbox{\rm a.s.}\leqno(2.11)
$$

\proclaim Fact 2.5. {\rm (Strassen [18])}
If $f\in {\cal S}$ defined by (1.5), then for any partition
$x_0=0<x_1<\ldots<x_{k}<x_{k+1}=1$ we have
$$
\sum_{i=1}^{k+1}{(f(x_i)-f(x_{i-1}))^2\over x_i-x_{i-1}}\leq 1.
\leqno(2.12)
$$

\proclaim Fact 2.6. {\rm (Chung [6])}
$$
\liminf_{t\to\infty}\sqrt{\log\log t\over t}
\sup_{0\leq s\leq t}|W(s)|={\pi\over\sqrt{8}}, \qquad
\hbox{\rm a.s.}\leqno(2.13)
$$

Define $g(T):=\max\{s\leq T:\, W(s)=0\}$. A joint lower class
result for $g(T)$ and $M(T):=\sup_{0\leq s\leq T}|W(s)|$ reads as
follows.

\proclaim Fact 2.7. {\rm (Grill [15])}
{\sl Let $\beta(t), \, \gamma(t)$ be positive functions slowly varying at
infinity, such that $0<\beta(t)\leq 1$, $0<\gamma(t)\leq 1$, $\beta(t)$
is non-increasing, $\beta(t)\sqrt{t}\uparrow\infty$, $\gamma(t)$ is
monotone, $\gamma(t)t\uparrow\infty$, $\gamma(t)/\beta^2(t)$ is monotone.
Then}
$$
\p\left(M(T)\leq \beta(T)\sqrt{T},\, g(T)\leq \gamma(T)T \quad
{\rm i.o.}\right)=0\quad {\rm or}\quad 1
$$
{\sl according as $I(\beta,\gamma)<\infty$ or $=\infty$, where}
$$
I(\beta,\gamma)=\int_1^\infty
{1\over t\beta^2(t)}
\left(1+{\beta^2(t)\over\gamma(t)}\right)^{-1/2}
\exp\left(-{(4-3\gamma(t))\pi^2\over 8\beta^2(t)}\right) \d t.
$$

Now define $d(T):=\min\{s\geq T:\, W(s)=0\}$. Since $\{ d(T)> t\}=
\{ g(t) < T\}$, we deduce from Fact 2.7 the following estimate on
$d(T)$ when $T \to\infty$.

\proclaim Fact 2.8.
{\sl With probability 1}
$$
d(T)=O(T(\log T)^3),\qquad T\to\infty.
$$

\bigskip
\noindent {\section 3. Probability estimates}

\bigskip
\proclaim Lemma 3.1. For $T\geq 1$, $\delta, z >0$ we have
$$
\eqalign{
& \p\left(\sup_{0\le t\le T-1}\sup_{0\le s\le 1}
|Y(t+s)-Y(t)|>z\right)
\cr &\qquad\le
c_6\left(\sqrt{T}\exp\left(-{z^2\over 8(1+\delta)}\right)+
T\exp\left(-{z^2\over 2(1+\delta)}\right)\right)
\cr}
\leqno(3.1)
$$
with some positive constant $c_6=c_6(\delta)$.

\bigskip

For the proof see Cs\'aki et al. [7], Lemma 2.8.

\bigskip

\proclaim
Lemma 3.2. For $T>1$, $0<\delta<1/2$, $z>1$ we have
$$
\eqalign{
 &\p\left( \, \sup_{0\le t\le T-1}(Y(t+1)-Y(t)) \ge z\right)
\cr & \qquad \ge \min\left({1\over 2},\, {c\sqrt{T-1}\over z}
\exp\left(-{z^2\over 8(1-\delta)}\right)\right)-\exp\left(-z^2\right)
\cr }
\leqno(3.2)
$$
with some positive constant $c_7=c_7(\delta)>0$.

\noindent {\bf Proof.}
Let us construct an increasing sequence of stopping times by $\eta_0:=0$
and
$$
\eta_{k+1}: = \inf\{ t>\eta_k+1: \, W(t)=0\}, \qquad k=0,1,2,\dots
$$

\noindent Let
$$
\eqalign{
 &\nu_t:=\min\{i\ge 1: \, \eta_i>t\}
\cr & Z_i:=Y(\eta_{i-1}+1)-Y(\eta_{i-1}),\qquad i=1,2,\dots \cr }
$$

Then $(Z_i,\, \eta_i-\eta_{i-1})_{i\ge 1}$ are i.i.d. random
vectors with
$$
\eta_i-\eta_{i-1}\; \buildrel{\rm law}\over {=}\; 1+\tau^2,\qquad
Z_i\; \buildrel{\rm law}\over {=}\; Y(1),
$$
where $\tau$ has Cauchy distribution.  Clearly, for $t>0$,
$$
\sup_{0\le s\le t}(Y(s+1)-Y(s))\ge \max_{1\le i\le \nu_t} Z_i =
\overline Z_{\nu_t},
$$
with  $\overline Z_k:=\max_{1\le i\le k}Z_i$. First consider the
Laplace transform $(\lambda>0)$:
$$
\eqalign{ & \lambda\int_0^\infty e^{-\lambda u}\p\left(\overline
Z_{\nu_u}<z\right)\, \d u
    \cr &=\lambda \sum_{k=1}^\infty\e
\int_0^\infty e^{-\lambda u} 1_{\{\eta_{k-1}\le u<\eta_k\}}
1_{\{\overline Z_k<z\}}\, \d u
        \cr &=\sum_{k=1}^\infty\e\left( \Big[e^{-\lambda\eta_{k-1}}-e^{-\lambda\eta_k}\Big] 1_{\{\overline Z_k<z\}}\right)
                \cr &=\sum_{k=1}^\infty \left(\e \Big[1_{\{\overline
Z_k<z\}}e^{-\lambda\eta_{k-1}}\Big] -\e \Big[1_{\{\overline
Z_k<z\}} e^{-\lambda\eta_k}\Big]\right)
             \cr &=\sum_{k=1}^\infty\left(\e \Big[1_{\{\overline Z_{k-1}<z\}}
e^{-\lambda\eta_{k-1}}\Big] -\e \Big[1_{\{\overline Z_{k-1}<z,\,
Z_k\ge z\}} e^{-\lambda\eta_{k-1}}\Big] -\e \Big[1_{\{\overline
Z_k<z\}} e^{-\lambda\eta_k}\Big]\right)
             \cr &=1-\sum_{k=1}^\infty\e \Big[ 1_{\{\overline Z_{k-1}<z,\, Z_k\ge
             z\}} \, e^{-\lambda\eta_{k-1}} \Big]
             \cr &=1-\sum_{k=1}^\infty\e \Big[1_{\{\overline Z_{k-1}<z\} } e^{-\lambda\eta_{k-1}}\Big]\,
\p(Y(1)\ge z)
            \cr &=1-\sum_{k=1}^\infty \left(\e \Big[
1_{\{Z_1<z\}} e^{-\lambda\eta_1}\Big]\, \right)^{k-1}\, \p (Y(1)
\ge z) \cr &=1-{\p(Y(1)\ge z)\over 1-\e \Big[ 1_{\{Z_1<z\}}
e^{-\lambda\eta_1}\Big]}, \cr}
$$
i.e.,
$$
\lambda\int_0^\infty e^{-\lambda u}\p\left(\overline Z_{\nu_u}\ge
z\right)\, \d u= {\p(Y(1)\ge z)\over 1-\e \Big[ 1_{\{Z_1<z\}}
e^{-\lambda\eta_1}\Big]}. \leqno(3.3)
$$

But (recalling that $Z_1=Y(1)$)
$$
1-\e \Big[ 1_{\{Z_1<z\}} e^{-\lambda\eta_1}\Big] \le 1-\e
(e^{-\lambda\eta_1})+\p (Y(1)\ge z)
$$
and (cf. [14], 3.466/1)
$$
1-\e e^{-\lambda\eta_1}=1-{1\over \pi}\int_{-\infty}^\infty
{e^{-\lambda(1+x^2)}\over 1+x^2}\, \d x={2\over\sqrt{\pi}}
\int_0^{\sqrt{\lambda}}e^{-x^2}\, \d x\le 2\sqrt{\lambda},
$$
hence
$$
\lambda\int_0^\infty e^{-\lambda u}\p\left(\overline
Z_{\nu_u}\ge z\right)\, \d u\ge {\p(Y(1)\ge z)\over
2\sqrt{\lambda}+\p(Y(1)\ge z)}.
$$

On the other hand, for any  $u_0>0$ we have
$$
\eqalign{  \lambda\int_0^\infty e^{-\lambda u}\p\left(\overline
Z_{\nu_u}\ge z\right)\, \d u &= \lambda\int_0^{u_0} e^{-\lambda
u}\p\left(\overline Z_{\nu_u}\ge z\right)\, \d u+
\lambda\int_{u_0}^\infty e^{-\lambda u}\p\left(\overline
Z_{\nu_u}\ge z\right)\, \d u \cr &\le \p\left(\overline
Z_{\nu_{u_0}}\ge z\right)+e^{-\lambda u_0}. \cr}
$$

\noindent It turns out that
$$
\p\left(\overline Z_{\nu_{u_0}}\ge z\right)\ge
{\p(Y(1)\ge z)\over
2\sqrt{\lambda}+\p(Y(1)\ge z)}-e^{-\lambda u_0}
\ge \min\left({1\over 2},\, {\p(Y(1)\ge z)\over
4\sqrt{\lambda}}\right)-e^{-\lambda u_0},
$$
where the inequality
$$
{x\over y+x}\ge \min\left({1\over2},\, {x\over 2y}\right),\qquad x> 0,\,
y>0
$$
was used. Choosing $u_0=T-1$, $\lambda=z^2/u_0$, and applying
(2.5) of Fact 2.1, we finally get
$$
\eqalign{
&\p\left( \, \sup_{0\le t\le T-1}(Y(t+1)-Y(t)) \ge
z\right)
\cr & \quad
\ge \min\left({1\over2},\, {c_8(\delta)\sqrt{T-1}\over
z}\exp\left(-{z^2\over
8(1-\delta)}\right)\right)
-\exp\left(-z^2\right).
\cr}
\leqno(3.4)
$$
This proves Lemma 3.2. \qed

\bigskip
\proclaim Lemma 3.3. For $T\ge 2$, $0\le \kappa< 1$  and  $\delta,
z>0$ we have
$$
\p\left(\sup_{0\le t\le T-1}(Y(t+1)-Y(t))<z\right)
\le {5\over T^{\kappa/2}}+\exp\left(-c_9T^{(1-\kappa)/2}
e^{-(1+\delta)z^2/8}\right)
\leqno(3.6)
$$
with some positive constant $c_9=c_9(\delta)$.

\medskip
See Cs\'aki et al. [7], Lemma 3.1.

\medskip
\proclaim Lemma 3.4. For $T>1$, $0<z\le 1/2$ we have
$$
\p\left(\sup_{0\le t\le T-1}\sup_{0\le s\le 1} |Y(t+s)-Y(t)|<z\right)
\ge {c_{10}\over\sqrt{T}}\exp\left(-{c_{11}\over z^2}\right)
$$
with some positive constants $c_{10},c_{11}$.

\medskip\noindent
{\bf Proof.} Define the events
$$
A:=\left\{\sup_{0\le s\le 1}|Y(s)|<{z\over 4},\, W(1)\ge {4\over z},
\inf_{1\le u\le T} W(u)\ge {2\over z}\right\}
$$
and
$$
\widetilde A:=\left\{\sup_{0\le t\le T-1}\sup_{0\le s\le 1}
|Y(t+s)-Y(t)|<z\right\}.
$$
Then $A\subset \widetilde A$, since if $A$ occurs and $t<1$,
$t+s\le 1$, then
$$
|Y(t+s)-Y(t)|\le 2\sup_{0\le s\le 1}|Y(s)|\le {z\over 2}<z.
$$

If $A$ occurs and $t<1$, $s\le 1$, $1<t+s\le T$, then
$$
|Y(t+s)-Y(t)|\le Y(t+s)-Y(1)+|Y(t)-Y(1)|\le \int_1^{t+s}{\d u\over W(u)}
+{z\over 2}< z.
$$

Moreover, if $A$ occurs and $1\le t$, $s\le 1$, $t+s\le T$, then
$$
|Y(t+s)-Y(t)|=\int_t^{t+s} {\d u\over W(u)}\le {z\over 2}<z.
$$

Hence $A\subset \widetilde A$ as claimed. But by the Markov
property of $W$,
$$
\p(A)=\int_{4/z}^\infty \p \left(\sup_{0\le s\le 1}|Y(s)|<{z\over 4}
\, \Big|\, W(1)=x\right)\p \left(\inf_{1\le u\le T} W(u)\ge {2\over z}
\, \Big|\, W(1)=x\right)\varphi(x)\, \d x,
\leqno(3.8)
$$
where $\varphi$ denotes the standard normal density function.

Using reflection principle and $x\ge 4/z$, $z\le 1/2$, we get
$$
\eqalign{
& \p\left(\inf_{1\le u\le T}W(u)\ge {2\over z}\, \Big|\, W(1)=x\right)=
2\Phi\left({x-2/z\over\sqrt{T-1}}\right)-1
\cr & \ge 2\Phi\left({2\over z\sqrt{T-1}}\right)-1\ge
2\Phi\left({4\over\sqrt{T}}\right)-1\ge {c_{12}\over\sqrt{T}},\cr}
\leqno(3.9)
$$
with some constant $c>0$, where $\Phi(\cdot)$ is the standard
normal distribution function. Hence
$$
\p(\widetilde A)\ge\p(A)\geq
{c_{12}\over\sqrt{T}}\p\left(\sup_{0\le s\le 1}|Y(s)|\le {z\over
4},\, W(1)\ge {4\over z}\right).\leqno(3.10)
$$
To get a lower bound of the probability on the right-hand side,
define $g$, $(m(v),\, 0\le v\le 1)$, $(B(u),\, 0\le u\le 1)$ by
(2.1), (2.2) and (2.3), respectively. Recall (see Fact 2.1 ) that
these three objects are independent, $g$ has arc sine
distribution, $m$ is a Brownian meander and $B$ is a Brownian
bridge. Moreover, $(g, m , B)$ are independent of ${\tt
sgn}(W(1))$ which is a Bernoulli variable. Observe that
$$ \eqalign{
\sup_{0\le s\le g} |Y(s)| &=\sqrt{g}\sup_{0\le s\le
1}\left|\int_0^s {\d u\over B(u)}\right|,
     \cr
\sup_{g\le s\le 1} |Y(s)| &= |Y(1)-Y(g)| =\sqrt{1-g}\int_0^1{\d
v\over m(v)},
   \cr |W(1)|&=\sqrt{1-g}\, m(1).\cr}
$$

Then
$$
\eqalign{
&\p\left(\sup_{0\le s\le 1}|Y(s)|\leq{z\over 4},\, W(1)\ge {4\over
z}\right)
   \cr & \ge \p\left(\sup_{0\le s\le g}|Y(s)|\leq{z\over 8},\, Y(1)-Y(g)\leq{z\over
8},\, W(1)\ge {4\over z}\right)
    \cr &\ge \p\left(\sqrt{g}\sup_{0\le s\le 1}\left|\int_0^s {\d u\over
B(u)}\right|\leq{z\over 8},\, \sqrt{1-g}\int_0^1{\d v\over
m(v)}\leq{z\over 8},\, \sqrt{1-g}\, m(1)\ge {4\over z}, \, W(1)
>0, \,  g<z^2\right)
    \cr &\ge \p\left(\sup_{0\le s\le 1}\left|\int_0^s {\d u\over
B(u)}\right|\leq{1\over 8},\, \int_0^1{\d v\over m(v)}\leq{z\over
8},\, m(1)\ge {4\over z\sqrt{1-z^2}},\,W(1)
>0, \, g<z^2\right)
    \cr &=\p\left(\sup_{0\le s\le 1}\left|\int_0^s {\d u\over
B(u)}\right|\leq{1\over 8}\right) \p\left(\int_0^1{\d v\over
m(v)}\leq{z\over 8},\, m(1)\ge {4\over z\sqrt{1-z^2}}\right)
\p(W(1) >0)\p(g<z^2)
    \cr &\ge c_{13}z\p\left(\int_0^1{\d v\over m(v)}\leq{z\over
8},\, m(1)\ge {4\over z\sqrt{1-z^2}}\right)
    \cr &=c_{13}z\int_{4/(z\sqrt{1-z^2})}^\infty \p\left(\int_0^1{\d v\over
m(v)}\leq{z\over
8}\, \Big|\, m(1)=x\right)\p(m(1)\in \d x).
\cr}
$$

It follows from Facts 2.1 and 2.2 that for $x>0$, $z>0$
$$
\p\left(\int_0^1{\d v\over m(v)}\leq{z\over
8}\, \Big|\, m(1)=x\right)\ge\p\left(\int_0^1{\d v\over m(v)}\leq{z\over
8}\, \Big|\, m(1)=0\right)\ge {c_{14}\over z^3}\exp\left(-{c_{15}\over
z^2}\right)
\leqno(3.11)
$$
and
$$
\p\left(m(1)>{4\over z\sqrt{1-z^2}}\right)=\exp\left(-{8\over z^2(1-z^2)}
\right).\leqno(3.12)
$$
Putting (3.10), (3.11), (3.12) together, we get (3.7). \qed

\proclaim
Lemma 3.5. For $T>1$, $0<z\leq 1/2$, $0<\delta\leq 1/2$ we have
$$
\eqalign{
& \p\left(\inf_{0\le t\le T-1}\sup_{0\le s\le 1}|Y(t+s)-Y(t)|<z\right)
\cr &\le c_{16}\left(\exp\left(-{(1-\delta)^2\over
2(1+\delta)^2z^2T}\right)
+\exp\left(-{c_5\delta\over 4(1+\delta)^2 z^2}\right)
+\exp\left({c_{17}\over z^2}-{c_{18}z^2\over
T}e^{c_{19}/z^2}\right)\right)\cr}
\leqno(3.13)
$$
with some positive constants $c_{16}$, $c_{17}=c_{17}(\delta)$,
$c_{18}=c_{18}(\delta)$, $c_{19}=c_{19}(\delta)$.

\bigskip

{\noindent\bf Proof.} Consider a positive integer $N$ to be given
later, $h=(T-1)/N$, $t_k=kh$, $k=0,1,2,\ldots, N$.
Then for $0<\delta\leq 1/2$ we have
$$
\leqalignno{ &\p\left(\inf_{0\le t\le T-1}\sup_{0\le s\le
1}|Y(t+s)-Y(t)|<z\right) \cr & \le \p\left( \inf_{0\le k\le N}
\sup_{0\le s\le 1} |Y(t_k+s)-Y(t_k)| \le (1+\delta)z\right)
+\p\left( \sup_{0\le t\le T-1} \sup_{0\le s\le h } | Y(t+s) -
Y(t)| > \delta z\right)
    \cr & =: P_1 + P_2 .
     \cr}$$
By scaling and Lemma 3.1
$$ \leqalignno{
P_2 &=\p\left(\sup_{0\le t\le (T-1)/h}\sup_{0\le s\le
1}|Y(t+s)-Y(t)|> {\delta z\over \sqrt{h}}\right)
     \cr&\le  c_6\left( \sqrt{ {T-1\over h}+1} \, \exp\left( -
{\delta^2 z^2 \over 8 h(1+\delta)}\right) + \left({T-1\over h}+1\right)
\exp\left(-{\delta^2 z^2\over 2h(1+\delta)}\right)\right)
     \cr&
\le 2c_6(N+1)\exp\left(-{\delta^2 z^2\over 8h(1+\delta)}\right).
    \cr}
$$

To bound $P_1$, we denote by $d(t):= \inf\{ s\ge t:
W(s)=0\}$ the first zero of $W$ after $t$.  Consider those $k$
for which $ \sup_{0\le s\le 1} | Y(t_k+s) - Y(t_k)| \le
(1+\delta)z$. If, moreover, $ d(t_k) \ge t_k+1 - \delta  $,
which means that the Brownian motion $W$ does not change sign over
$[t_k, t_k+1-\delta)$, then
$$ (1+\delta)z  \ge  |Y(t_k+1
- \delta ) - Y(t_k)| = \int_0^{1-\delta} { \d s \over | W(t_k+s)|}
\ge {1-\delta \over \sup_{0\le s\le T} |W(s)|},
$$
and it follows that
$$
\leqalignno{& P_1  \le \p\left(\sup_{0\le s\le T} |W(s)| >
{(1-\delta)
\over z(1+\delta)} \right)  \cr & + \p\left( \exists
k\le N : \sup_{0\le s\le 1} | Y(t_k+s) - Y(t_k)| \le (1+\delta)z;
d(t_k) < t_k+1 - \delta \right)
      \cr & \le 4 \exp\left(-{(1-\delta)^2\over 2(1+\delta)^2 z^2T}
      \right)
\cr & + \sum_{k=0}^{N}  \p\left(\sup_{0\le s\le 1} | Y(t_k+s)
- Y(t_k)|
      \le (1+\delta)z; d(t_k) < t_k+1 - \delta \right).
       \cr}
$$

Let $\widehat W(s) = W(d(t_k) +s)$ for $s\ge0$ and
$\widehat Y(s)$ be the associated principal values.  Observe that
on $\{\sup_{0\le s\le 1} | Y(t_k+s) - Y(t_k)| \le (1+\delta)z; d(t_k) <
t_k+1 - \delta \}$, we have $ \sup_{ 0\le u \le
\delta } | \widehat Y(u) + (Y(d(t_k)) - Y(t_k)) | <(1+\delta)z,$
and $ |Y(d(t_k)) - Y(t_k) | \le (1+\delta)z$
which implies that
$$
\sup_{ 0\le u \le \delta } | \widehat Y(u) |
<  2 (1+\delta)z.
$$
By scaling and Fact 2.3 we have
$$
\p\left( \sup_{ 0\le u \le \delta } | \widehat Y(u) | <  2
(1+\delta)z \right) \le c_4\exp\left( - { c_5 \delta \over
4(1+\delta)^2 z^2}\right).
$$
\noindent Therefore, we obtain:
$$
P_1 \le 4 \exp\left(-{(1-\delta)^2\over 2(1+\delta)^2z^2T}\right) +
c_4(N+1)\exp\left( - { c_5 \delta \over 4(1+\delta)^2z^2}\right).
$$
Hence  $$ P_1+P_2 \le 4\exp\left(-{(1-\delta)^2\over 2(1+\delta)^2
z^2  T}\right) + c_4(N+1)  \exp\left( - { c_5 \delta \over
4(1+\delta)^2z^2}\right)
$$
$$ + 2c_6(N+1)\exp\left(-{\delta^2 z^2\over
8h(1+\delta)}\right).
$$
By taking $N= [e^{c_5 \delta /(4(1+\delta)^2z^2)}]+1$, we get
$$
\eqalign{
& P_1+P_2
\cr &\le c_{16}\left(\exp\left(-{(1-\delta)^2\over
2(1+\delta)^2z^2T}\right)
+\exp\left(-{c_5\delta\over 4(1+\delta)^2 z^2}\right)
+\exp\left({c_{17}\over z^2}-{c_{18}z^2\over
T}e^{c_{19}/z^2}\right)\right)\cr}
$$
with relevant constants $c_{16}$, $c_{17}$, $c_{18}$, $c_{19}$,
proving (3.13).
\qed

\bigskip
\bigskip
\noindent {\section 4. Proof of Theorem 1.1(i)}

\bigskip
The upper estimation, i.e.
$$
\limsup_{T\to\infty} {\sup_{0\le t\le T-a_T}\; \sup_{0\le s\le a_T} \;
|Y(t+s)-Y(t)|\over \sqrt{8a_T \left(\log \sqrt {T/a_T}+\log\log T\right)}}
\le 1, \qquad\hbox{\rm a.s.}
\leqno(4.1)
$$
follows easily from Wen's Theorem E.

Now we prove the lower bound, i.e.
$$
\limsup_{T\to\infty} {\sup_{0\le t\le T-a_T}\; \sup_{0\le s\le a_T} \;
|Y(t+s)-Y(t)|\over \sqrt{8a_T
\left(\log \sqrt {T/a_T}+\log\log T\right)}} \ge 1, \qquad
\hbox{\rm a.s.}
\leqno(4.2)
$$

In the case when $a_T=T$, (4.2) follows from the law of the
iterated logarithm (1.3) of Theorem A. Now we assume that $a_T/T\le
\rho<1$, with some constant $\rho$ for all $T>0$.

By scaling, (3.2) of Lemma 3.2 is equivalent to
$$
\eqalign{
 &\p\left( \, \sup_{0\le t\le T-a}(Y(t+a)-Y(t)) \ge z\sqrt{a}\right)
\cr & \qquad \ge \min\left({1\over 2},\, {c_7\sqrt{T/a-1}\over z}
\exp\left(-{z^2\over 8(1-\delta)}\right)\right)-\exp\left(-z^2\right)
\cr }
\leqno(4.3)
$$
for $0<a<T$, $0<\delta<1/2$, $z>1$.

Define the sequences
$$
t_k:=e^{7k\log k},\qquad k=1,2,\dots
\leqno(4.4)
$$
and $\theta_0:=0$,
$$
\theta_k:=\inf\{t>T_k:\, W(t)=0\},\qquad k=1,2,\dots,
\leqno(4.5)
$$
where $T_k:=\theta_{k-1}+t_k$. For $0<\delta<\min(1/2,1-\rho)$ define the
events
$$
A_k:=\left\{\sup_{0\le t\le t_k(1-\delta)-a_{t_k}}
(Y(\theta_{k-1}+t+a_{t_k})-Y(\theta_{k-1}+t))\ge
(1-\delta)\beta_k\right\},\quad k=1,2,\dots
$$
with
$$
\beta_k:=\sqrt{8a_{t_k}\left(\log\sqrt{t_k\over a_{t_k}}+\log\log
t_k\right)}.
$$

Applying (4.3) with $T=t_k(1-\delta)$, $a=a_{t_k}$,
$z=(1-\delta)\sqrt{8(\log\sqrt{t_k/a_{t_k}}+\log\log t_k)}$, we have
for $k$ large
$$
\eqalign{
&\p(A_k)=\p\left(\sup_{0\le t\le t_k(1-\delta)-a_{t_k}}
(Y(t+a_{t_k})-Y(t))\ge (1-\delta)\beta_k\right)
\cr & \quad
\ge \min\left({1\over 2},\, {b_k\over (\log
t_k)^{1-\delta}}\right)-{1\over (\log t_k)^{8(1-\delta)^2}}
\cr}
$$
with
$$
b_k={c_7\sqrt{t_k(1-\delta)/a_{t_k}-1}\over
(t_k/a_{t_k})^{(1-\delta)/2}\sqrt{\log\sqrt{t_k/a_{t_k}}+\log\log t_k}}
\ge {c_{20}\over \sqrt{\log k}}.
$$

Hence $\sum_k\p(A_k)=\infty$ and since $A_k$ are independent,
Borel-Cantelli lemma yields
$$
\p(A_k\, {\rm i.o.})=1.
$$
It follows that
$$
\limsup_{k\to\infty}
{\sup_{0\le t\le t_k(1-\delta)-a_{t_k}}
(Y(\theta_{k-1}+t+a_{t_k})-Y(\theta_{k-1}+t))\over
\sqrt{8a_{t_k}\left(\log\sqrt{t_k\over a_{t_k}}+\log\log
t_k\right)}}\ge 1-\delta,\qquad
\hbox{\rm a.s.}
\leqno(4.6)
$$

It can be seen (cf. [9]) that we have almost surely for large enough $k$
$$
t_k\le T_k\le t_k\left(1+{1\over k}\right),
$$
consequently
$$
\lim_{k\to\infty}{t_k\over T_k}=1,\qquad
\hbox{\rm a.s.}
\leqno(4.7)
$$
Since by our assumptions
$$
{t_k\over T_k}\le {a_{t_k}\over a_{T_k}}\le 1,
$$
we have also
$$
\lim_{k\to\infty}{a_{t_k}\over a_{T_k}}=1,\qquad
\hbox{\rm a.s.}
\leqno(4.8)
$$
On the other hand, for any $\delta>0$ small enough we have almost
surely for large $k$
$$
a_{T_k}\le (1+\delta)a_{t_k}\le t_k\delta+a_{t_k},
$$
thus
$$
T_k-a_{T_k}\ge T_k-t_k\delta-a_{t_k},
$$
consequently
$$
\eqalign{
& \sup_{0\le t\le T_k-a_{T_k}}\sup_{0\le s \le a_{T_k}}
|Y(t+s)-Y(t)|
\cr &\ge \sup_{0\le t\le t_k(1-\delta)-a_{t_k}}
(Y(\theta_{k-1}+t+a_{t_k})-Y(\theta_{k-1}+t)),
\cr}
\leqno(4.9)
$$
hence we have also
$$
\limsup_{k\to\infty}
{\sup_{0\le t\le T_k-a_{T_k}}\sup_{0\le s \le a_{T_k}}
|Y(t+s)-Y(t)|
\over
\sqrt{8a_{t_k}\left(\log\sqrt{t_k\over a_{t_k}}+\log\log
t_k\right)}}\ge 1-\delta,\qquad
\hbox{\rm a.s.}
\leqno(4.10)
$$
and since $\delta>0$ can be arbitrary small, (4.2) follows by combining
(4.7), (4.8), (4.9) and (4.10). \qed

\bigskip
\bigskip
\noindent {\section 5. Proof of Theorem 1.1(ii)}

\bigskip
First assume that
$$
a_T>{T\over (\log T)^\alpha} \qquad\hbox{\rm for some}\quad \alpha<2.
\leqno(5.1)
$$
By Theorem C,  $$ \eqalign{ & \liminf_{T\to\infty}\sqrt{\log\log
T\over a_T} \sup_{0\le t\le T-a_T}\sup_{0\le s\le
a_T}|Y(t+s)-Y(t)| \cr &\ge\liminf_{T\to\infty}\sqrt{\log\log
a_T\over a_T} \sup_{0\le s\le a_T}|Y(s)|\ge K_1, \qquad \hbox{\rm
a.s.},\cr} \leqno(5.2)
$$
proving the lower bound in (1.12).

To get an upper bound, note that by scaling, (3.7) of Lemma 3.4 is
equivalent to
$$
\p\left(\sup_{0\le t\le T-a}\sup_{0\le s\le a}
|Y(s+t)-Y(t)|<z\sqrt{a}\right)
\ge c_{10}\sqrt{a\over T}\exp\left(-{c_{11}\over z^2}\right)
\leqno(5.3)
$$
for $T\ge a$, $0<z\le 1/2$.

Let $t_k$ and $\theta_k$ be defined by (4.4) and (4.5), resp.,  as in the
proof of Theorem 1.1(i) and for any
$\varepsilon>0$ and for $\delta>0$ such that $\alpha/2+c_{11}/\delta^2<1$,
define the events
$$
E_k:=\left\{\sup_{0\le t\le
(1+\varepsilon)t_k-a_{t_k(1+\varepsilon)}} \sup_{0\le s\le
a_{t_k(1+\varepsilon)}} |Y(\theta_{k-1}+t+s)-Y(\theta_{k-1}+t)|
\le\delta\sqrt{a_{t_k}\over\log\log t_k}\right\}.
$$   Then putting $T=(1+\varepsilon)t_k$, $a=a_{(1+\varepsilon)t_k}$,
$z=\delta/\sqrt{\log\log t_k}$, into (5.3), we get $$
\p(E_k)=\p\left(\sup_{0\le t\le
(1+\varepsilon)t_k-a_{t_k(1+\varepsilon)}} \sup_{0\le s\le
a_{t_k(1+\varepsilon)}} |Y(t+s)-Y(t)|
\le\delta\sqrt{a_{t_k}\over\log\log t_k}\right)
$$
$$
\ge c_{10}\sqrt{a_{t_k}\over
t_k}\exp(-(c_{11}/\delta^2)\log\log((1+\varepsilon)t_k))\ge
{c_{10}\over (\log t_k)^{\alpha/2+c_{11}/\delta^2}},
$$
hence $\sum_k\p(E_k)=\infty$, and since $E_k$ are independent, we
have $\p(E_k\, {\rm i.o.})=1$, i.e.
  $$
\liminf_{k\to\infty}\sqrt{\log\log t_k\over a_{t_k}}\sup_{0\le t\le
(1+\varepsilon)t_k-a_{t_k(1+\varepsilon)}}
\sup_{0\le s\le a_{t_k(1+\varepsilon)}}
|Y(\theta_{k-1}+t+s)-Y(\theta_{k-1}+t)|
\le\delta,\quad{\rm a.s.}
\leqno(5.4)
$$
for any $\varepsilon$. Put, as before, $T_k=\theta_{k-1}+t_k$. For
large enough $k$ by (4.7) and (4.8) we have $ a_{T_k}\le
(1+\varepsilon)a_{t_k}, \,{\rm a.s.} $ and $
T_k-a_{T_k}\le\theta_{k-1}+(1+\varepsilon)t_k-(1+\varepsilon)a_{t_k},
\,{\rm a.s.} $ Thus given any $\varepsilon>0$, we have for large
$k$  $$ \eqalign{ & \sup_{0\le t\le T_k-a_{T_k}}\sup_{0\le s\le
a_{T_k}} |Y(t+s)-Y(t)|
     \cr & \le 2\sup_{0\le t\le \theta_{k-1}}|Y(t)| + \sup_{0\le t\le
(1+\varepsilon)t_k-a_{t_k(1+\varepsilon)}} \sup_{0\le s\le
a_{t_k(1+\varepsilon)}}
|Y(\theta_{k-1}+t+s)-Y(\theta_{k-1}+t)|.\cr} \leqno(5.5)
$$

By Theorem A, Fact 2.8, (4.7), (5.1) and simple calculation,
$$
\eqalign{
& \sup_{0\le t\le \theta_{k-1}}|Y(t)|=
O(\theta_{k-1}\log\log\theta_{k-1})^{1/2}
\cr & =O(t_{k-1}(\log t_{k-1})^3\log\log t_{k-1})^{1/2}
=o\left({a_{t_k}\over\log\log t_k}\right)^{1/2},\qquad{\rm a.s.}\cr}
\leqno(5.6)
$$
as $k\to\infty$.
  Assembling
(5.4), (5.5) and (5.6), we get
$$
\liminf_{k\to\infty}\sqrt{\log\log t_k\over a_{t_k}}
\sup_{0\le t\le T_k-a_{T_k}}\sup_{0\le s\le a_{T_k}}|Y(t+s)-Y(t)|
$$
$$
=\liminf_{k\to\infty}\sqrt{\log\log T_k\over a_{T_k}}
\sup_{0\le t\le T_k-a_{T_k}}\sup_{0\le s\le a_{T_k}}|Y(t+s)-Y(t)|
\le\delta,\qquad{\rm a.s.}
$$
which together with (5.2) yields (1.12).

Now assume that
$$
a_T\le {T\over(\log T)^\alpha}\qquad\hbox{\rm for some}\quad \alpha>2.
\leqno(5.7)
$$

By Theorem 1.1(i),
$$
\eqalign{
& \liminf_{T\to\infty}
{\sup_{0\le t\le T-a_T}\sup_{0\le s\le a_T}|Y(t+s)-Y(t)|
\over\sqrt{a_T\log(T/a_T)}}
\cr &\le\limsup_{T\to\infty}
{\sup_{0\le t\le T-a_T}\sup_{0\le s\le a_T}|Y(t+s)-Y(t)|
\over\sqrt{a_T\log(T/a_T)}}
\cr &\le\limsup_{T\to\infty}
{\sup_{0\le t\le T-a_T}\sup_{0\le s\le a_T}|Y(t+s)-Y(t)|
\over\sqrt{{2\alpha a_T\over \alpha+2}\left(\log\sqrt{T/a_T}+\log\log
T\right)}}\le2\sqrt{\alpha+2\over\alpha},\cr}
\leqno(5.8)
$$
i.e., an upper bound in (1.13) follows.

To get a lower bound under (5.7), observe that by scaling, (3.6) of Lemma
3.3 is equivalent to
$$
\p\left(\sup_{0\le t\le T-a}(Y(t+a)-Y(t))<z\sqrt{a}\right)
\le 5\left({a\over T}\right)^{\kappa/2}+\exp\left(-c_9\left({T\over
a}\right)^{(1-\kappa)/2}
e^{-(1+\delta)z^2/8}\right)
$$
for $a\le T$, $0\le\kappa<1$, $0<\delta$, $0<z$. Using (5.7) we get
further
$$
\eqalign{
& \p\left(\sup_{0\le t\le T-a}(Y(t+a)-Y(t))<z\sqrt{a}\right)
\cr &\le {5\over (\log T)^{\alpha\kappa/2}}+
\exp\left(-c_9\left(\log T)^{\alpha(1-\kappa)/2}\right)
e^{-(1+\delta)z^2/8}\right).\cr}
\leqno(5.9)
$$
In the case when (1.7) holds, (1.13) was proved in [7]. In other cases the
proof is similar. Let $T_k=e^k$ and define the events
$$
F_k=\left\{\sup_{0\leq t\leq T_k-a_{T_k}} (Y(t+a_{T_k})-Y(t))\leq
C_1\sqrt{a_{T_k}\log{T_k\over a_{T_k}}}\right\}
$$
with some constant $C_1$ to be given later. By (5.9)
$$
\p(F_k)\leq {5\over k^{\alpha\kappa/2}}+
\exp\left(-c_9k^{\alpha((1-\kappa)/2-(1+\delta)C_1^2/8)}\right).
$$
For given $\alpha>2$, choose small $\varepsilon>0$,
$\kappa=2/\alpha+\varepsilon$,
$$
C_1=2\sqrt{\alpha-2-2\varepsilon(1+\alpha)\over (1+\varepsilon)\alpha}.
$$
One can easily see that with these choices $\sum_k\p(F_k)<\infty$,
consequently
$$
\liminf_{k\to\infty}{\sup_{0\leq t\leq T_k-a_{T_k}} (Y(t+a_{T_k})-Y(t))
\over \sqrt{a_{T_k}\log{T_k\over a_{T_k}}}}
\geq C_1, \qquad{\rm a.s.},
$$
implying also
$$
\liminf_{k\to\infty}{\sup_{0\leq t\leq T_k-a_{T_k}}
\sup_{0\leq s\leq a_{T_k}} |Y(t+s)-Y(t)|
\over \sqrt{a_{T_k}\log{T_k\over a_{T_k}}}}
\geq 2\sqrt{\alpha-2\over\alpha}, \qquad{\rm a.s.},
$$
for $\varepsilon$ can be choosen arbitrary small.

Since $\sup_{0\leq t\leq T-a_T}\sup_{0\leq s\leq a_T}
|Y(t+s)-Y(t)|$ is increasing in $T$, we obtain  a lower bound in
(1.13). This together with the  0-1 law for Brownian motion
complete  the proof of Theorem 1.1(ii).\qed

\bigskip
\noindent {\section 6. Proof of Theorem 1.2(i)}

\bigskip\noindent
If $a_T=T$, then (1.14) is equivalent to Theorem C. Now assume that
$\rho:=\lim_{T\to\infty}a_T/T<1$.

First we prove the lower bound, i.e.
$$
\liminf_{T\to\infty}\; { \sqrt{ T \log\log T} \over a_T}
\inf_{0\le t\le T-a_T}\; \sup_{0\le s\le a_T} \; |Y(t+s)-Y(t)|\geq c,
\qquad \hbox{\rm a.s.}\leqno(6.1)
$$

By scaling, (3.13) of Lemma 3.5 is equivalent to
$$
\eqalign{
& \p\left(\inf_{0\le t\le T-a}\sup_{0\le s\le a}|Y(t+s)-Y(t)|<z\right)
\cr &\le c_{16}\left(\exp\left(-{a(1-\delta)^2\over
2(1+\delta)^2z^2T}\right)
+\exp\left(-{c_5\delta\over 4(1+\delta)^2 z^2}\right)
+\exp\left({c_{17}\over z^2}-{c_{18}az^2\over
T}e^{c_{19}/z^2}\right)\right)\cr}
\leqno(6.2)
$$
for $a<T$, $0<z\leq 1/2$, $0<\delta\leq 1/2$.

Define the events
$$
G_k=\left\{\inf_{0\le t\le T_{k+1}-a_{T_k}}\sup_{0\le s\le
a_{T_k}}|Y(t+s)-Y(t)|<z_k\right\}\quad k=1,2,\dots
$$
Let $T_k=e^k$ and put $T=T_{k+1}$, $a=a_{T_k}$,
$$
z=z_k=C_2\sqrt{a_{T_k}\over T_{k+1}\log\log T_{k+1}}
$$
into (6.2). The constant $C_2$ will be choosen later.
Denoting the terms on the right-hand side of (6.2) by $I_1$,
$I_2$, $I_3$, resp., we have
$$
\p(G_k)\leq c_{16}(I_1^{(k)}+I_2^{(k)}+I_3^{(k)}),
$$
where
$$
I_1^{(k)}=\exp\left(-{c_{21}\over C_2^2}\log\log T_{k+1}\right),
$$
$$
I_2^{(k)}=\exp\left(-{c_{22}T_k\over C_2^2 a_{T_k}}\log\log
T_{k+1}\right),
$$
$$
I_3^{(k)}=\exp\left({c_{23}T_k\log\log T_{k+1}\over C_2^2a_{T_k}}
-{c_{24}C_2^2a^2_{T_k}\over T_k^2\log\log T_{k+1}}
\left(\log T_{k+1}\right)^{c_{25}T_k\over C_2^2a_{T_k}}\right)
$$
with some constants $c_{21}=c_{21}(\delta)$, $c_{22}=c_{22}(\delta)$,
$c_{23}$, $c_{24}$, $c_{25}$.

One can see easily that for any choice of positive $C_2$ and for all
possible $a_T$ (satisfying our conditions) we have $\sum_k
I_3^{(k)}<\infty$. So we show that for appropriate choice of $C_2$ we have
also $\sum_kI_j^{(k)}<\infty$, $j=1,2$.

First consider the case $0<\rho>0$. Choosing a positive $\delta$
one can select $C_2<\min(\sqrt{c_{21}},\sqrt{c_{22}\over\rho})$
and it is easy to verify that $\sum_kI_j^{(k)}<\infty$, $j=1,2$,
hence also $\sum_k\p(G_k)<\infty$.

In the case $\rho=0$ choose $C_2<(1-\delta)/((1+\delta)\sqrt{2})$.
With this choice we have $\sum_k I_1^{(k)}<\infty$ for arbitrary
$\delta>0$. Since $\lim_{k\to\infty}(T_k/a_{T_k})=\infty$, we have
also $\sum_k I_2^{(k)}<\infty$ and $\sum_k\p(G_k)<\infty$.
Borell-Cantelli lemma and interpolation between $T_k$'s finish the
proof of (6.1). We have also verified that in the case $\rho=0$
one can choose $C_2=1/\sqrt{2}$, since $\delta$ can be choosen
arbitrary small.

Now we turn to the proof of the upper bound, i.e.
$$
\liminf_{T\to\infty}\; { \sqrt{ T \log\log T} \over a_T}
\inf_{0\le t\le T-a_T}\; \sup_{0\le s\le a_T} \; |Y(t+s)-Y(t)|\leq C_3,
\qquad \hbox{\rm a.s.}\leqno(6.3)
$$
with some constant $C_3$.

If $\rho>0$, then
$$
\inf_{0\le t\le T-a_T}\; \sup_{0\le s\le a_T} \; |Y(t+s)-Y(t)|
\leq \sup_{0\leq s\leq a_T}|Y(s)|\leq \sup_{0\leq s\leq T}|Y(s)|
$$
and hence (6.3) with some positive constant $C_3$ follows from Theorem C.

If $\rho=0$, then let for any $\varepsilon>0$
$$
\lambda_T:=\inf\{t:\, |W(t)|=\sup_{0\leq s\leq T(1-\varepsilon)}|W(s)|\}.
\leqno(6.4)
$$
According to the law of the iterated logarithm, with probability one there
exists a sequence $\{T_i,\, i\geq 1\}$ such that
$\lim_{i\to\infty}T_i=\infty$ and
$$
|W(\lambda_{T_i})|\geq \sqrt{2T_i(1-\varepsilon)\log\log T_i}.
\leqno(6.5)
$$
But Fact 2.4 implies that for $\varepsilon>0$
$$
|W(\lambda_{T_i})-W(s)|\leq\sqrt{2(1+\varepsilon)\varepsilon T_i\log\log
T_i}, \quad \lambda_{T_i}\leq s\leq\lambda_{T_i}+\varepsilon T_i,
\quad i\geq 1.\leqno(6.6)
$$
Now assume that $W(\lambda_{T_i})>0$. The case when $W(\lambda_{T_i})<0$
is similar. Then (6.5) and (6.6) imply
$$
W(s)\geq\left(\sqrt{1-\varepsilon}-\sqrt{\varepsilon(1+\varepsilon)}\right)
\sqrt{2T_i\log\log T_i}, \quad \lambda_{T_i}\leq s\leq
\lambda_{T_i}+\varepsilon T_i.\leqno(6.7)
$$
$\rho=0$ implies that $a_T\leq \varepsilon T$ for any $\varepsilon>0$ and
large enough $T$, hence we have from (6.7) for large $i$
$$
\sup_{0\le s\le a_{T_i}}(Y(\lambda_{T_i}+s)-Y(\lambda_{T_i}))=
Y(\lambda_{T_i}+a_{T_i})-Y(\lambda_{T_i})=
\int_{\lambda_{T_i}}^{\lambda_{T_i}+a_{T_i}}{\d s\over W(s)}
$$
$$
\leq {a_{T_i}\over
\left(\sqrt{1-\varepsilon}-\sqrt{\varepsilon(1+\varepsilon)}\right)
\sqrt{2T_i\log\log T_i}}.
$$
Since $\varepsilon>0$ is arbitrary, (6.3) follows with
$C_3=1/\sqrt{2}$. This completes the proof of Theorem 1.2(i). \qed

\bigskip
\bigskip
\noindent {\section 7. Proof of Theorem 1.2(ii)}

\bigskip\noindent
If $\rho=1$, then (1.15) is equivalent to (1.3) of Theorem A. So we may
assume that $0<\rho<1$.

First we prove the upper bound
$$
\limsup_{T\to\infty}
{\inf_{0\le t\le T-\rho T}\; \sup_{0\le s\le \rho T} \; |Y(t+s)-Y(t)|
\over\sqrt{8T\log\log T}} \leq\rho,
\qquad \hbox{\rm a.s.}\leqno(7.1)
$$  Let $k$ be the largest integer for which $k\rho< 1$ and put
$x_i=i\rho$, $i=0,1,\ldots,k$, $x_{k+1}=1$. It suffices to show
that if $f\in {\cal S}$ defined by (1.5), then $$ \min_{1\leq
i\leq k+1}|f(x_i)-f(x_{i-1})|\leq \rho.
$$
Assume on the contrary that
$$
|f(x_i)-f(x_{i-1})|>\rho,\qquad \forall i=1,2,\ldots,k+1.
$$
Then
$$
\sum_{i=1}^{k+1}{(f(x_i)-f(x_{i-1}))^2\over x_i-x_{i-1}}
>\sum_{i=1}^k {\rho^2\over\rho} +{\rho^2\over 1-k\rho}=
k\rho +{\rho^2\over 1-k\rho}\geq 1,
$$
contradicting (2.12) of Fact 2.5. This proves (7.1).

The lower bound
$$
\limsup_{T\to\infty}
{\inf_{0\le t\le T-\rho T}\; \sup_{0\le s\le \rho T} \; |Y(t+s)-Y(t)|
\over\sqrt{8T\log\log T}} \geq\rho,
\qquad \hbox{\rm a.s.}\leqno(7.2)
$$
follows from the fact that by Theorem B the function $f(x)=x,\, 0\leq
x\leq 1$ is a limit point of
$$
{Y(xt)\over\sqrt{8T\log\log T}}
$$
and for this function
$$
\min_{0\leq x\leq 1-\rho} |f(x+\rho)-f(x)|=\rho.
$$
This completes the proof of Theorem 1.2(iia). \qed

Now assume that
$$
\lim_{T\to\infty}{a_T(\log\log T)^2\over T}=0.\leqno(7.3)
$$
Define $\lambda_T$ as in (6.4). Then according to Chung's LIL (cf. Fact
2.6)
$$
|W(\lambda_T)|\geq {\pi\over\sqrt{8}}(1-\varepsilon)
\sqrt{T\over \log\log T}\leqno(7.4)
$$
for every $T$ sufficiently large. But according to Fact 2.4,
$$
\eqalign{
& \sup_{0\leq s\leq a_T}|W(\lambda_T+s)-W(\lambda_T)|
\cr &\leq\sqrt{(2+\varepsilon)a_T(\log(T/a_T)+\log\log T)}
\leq\sqrt{(2+\varepsilon)\varepsilon T\over\log\log T}.\cr}
$$
Assuming $W(\lambda_T)>0$, we get
$$
W(\lambda_T+s)\geq W(\lambda_T)-
\sqrt{(2+\varepsilon)\varepsilon T\over\log\log T}
\geq c\sqrt{T\over\log\log T}.
$$
Hence
$$
\inf_{0\leq t\leq T-a_T}\sup_{0\leq s\leq a_T}
|Y(t+s)-Y(t)|\leq Y(\lambda_T+a_T)-Y(\lambda_T)
$$
$$
=\int_0^{a_T}{\d s\over
W(\lambda_T+s)}\leq {a_T\over c}\sqrt{\log\log T\over T}
$$
for all large $T$.

The case when $W(\lambda_T)<0$ is similar. This shows the upper bound in
(1.16).

For the lower bound we use Fact 2.6: with probability one
$$
g_T\leq {T\over (\log\log T)^2},\quad
\max_{0\leq u\leq T} |W(u)|\leq
{\pi\over\sqrt{2}}\sqrt{T\over\log\log T}\quad {\rm i.o.}\leqno(7.5)
$$
According to Theorem 1.2(i) for every large $T$ we have for any
$\varepsilon>0$ and sufficiently large $T$
$$
\eqalign{
\cr &\inf_{0\leq t\leq T(\log\log T)^{-2}}\sup_{0\leq s\leq a_T}
|Y(t+s)-Y(t)|
\cr &\geq {(K_4-\varepsilon) a_T\over
\sqrt{\left({T\over(\log\log T)^2}+a_T\right)\log\log T}}
\leq {(K_4-\varepsilon)a_T\over\sqrt{(1+\varepsilon)T\log\log T}}.
\cr}
\leqno(7.6)
$$
On the other hand, if $T(\log\log T)^{-2}\leq t\leq T-a_T$, then by (7.5)
$$
|Y(t+a_T)-Y(t)|=\int_t^{t+a_T}{\d s\over |W(s)|}\geq
{a_T\sqrt{2\log\log T}\over\pi\sqrt{T}}.
$$

Combining (7.6) and (7.7) we get for $\varepsilon>0$ and all large $T$
$$
\inf_{0\leq t\leq T-a_T}\sup_{0\leq s\leq a_T}|Y(t+s)-Y(t)|
\geq\min\left({K_4-\varepsilon\over\sqrt{1+\varepsilon}},
{\sqrt{2}\over\pi}\right){a_T\sqrt{\log\log T}\over T}.
$$

This shows the lower bound in (1.16). The proof of Theorem
1.2(iib) is complete by applying the  0-1 law for Brownian
motion.\qed

\bigskip
\bigskip
\noindent {\section Acknowledgements}

The authors are indebted to Marc Yor for useful remarks.
Cooperation between the authors was supported by the joint
French--Hungarian Intergovernmental Grant "Balaton"  (grant no. F-39/00).

\bigskip
\bigskip

\noindent {\section References}

\baselineskip=12pt

\bigskip
\item{[1]} Ait Ouahra, M. and Eddahbi, M.: Th\'eor\`emes limites pour
certaines fonctionnelles associ\'ees aux processus stables sur l'espace
de H\"older. {\it Publ. Mat.} {\bf 45} (2001), 371--386.

\medskip
\item{[2]} Bertoin, J.: On the Hilbert transform of the local times
of a L\'evy process. {\it Bull. Sci. Math.} {\bf 119} (1995),
147--156.

\medskip
\item{[3]} Bertoin, J.: Cauchy's principal value of local times of L\'evy
processes with no negative jumps via continuous branching processes. {\it
Electronic J. Probab.} {\bf 2} (1997), Paper No. 6, 1--12.

\medskip
\item{[4]} Biane, P. and Yor, M.: Valeurs principales associ\'ees
aux temps locaux browniens. {\it Bull. Sci. Math.} {\bf 111}
(1987), 23--101.

\medskip
\item{[5]} Boufoussi, B., Eddahbi, M. and Kamont, A.: Sur la
d\'eriv\'ee fractionnaire du temps local brownien. {\it Probab. Math.
Statist.} {\bf 17} (1997), 311--319.

\medskip
\item{[6]} Chung, K.L.: On the maximumpartial sums of sequences of
independent random variables. {\it Trans. Amer. Math. Soc.} {\bf 64}
(1948), 205--233.

\medskip
\item{[7]} Cs\'aki, E., Cs\"org\H o, M. F\"oldes, A. and Shi, Z.:
Increment sizes of the principal value of Brownian local time. {\it
Probab. Th. Rel. Fields} {\bf 117} (2000), 515--531.

\medskip
\item{[8]} Cs\'aki, E., Cs\"org\H o, M. F\"oldes, A. and Shi, Z.:
Path properties of Cauchy's principal values related to local time.
{\it Studia Sci. Math. Hungar.} {\bf 38} (2001), 149--169.

\medskip
\item{[9]} Cs\'aki, E. and F\"oldes, A.: A note on the stability of the
local time of a Wiener process. {\it Stoch. Process. Appl.} {\bf 25 }
(1987), 203--213.

\medskip
\item{[10]} Cs\'aki, E., F\"oldes, A. and Shi, Z.: A joint functional law
for the Wiener process and principal value. {\it Studia Sci. Math.
Hungar.} {\bf 40} (2003), 213--241.

\medskip
\item{[11]} Cs\'aki, E., Shi, Z. and Yor, M.: Fractional Brownian motions
as "higher-order" fractional derivatives of Brownian local times. In: {\it
Limit Theorems in Probability and Statistics} (I. Berkes et al., eds.)
Vol. I, pp. 365--387. J\'anos Bolyai Mathematical Society, Budapest, 2002.

\medskip
\item{[12]} Cs\"org\H o, M. and R\'ev\'esz, P.: {\it Strong Approximations
in Probability and Statistics.} Academic Press, New York, 1981.

\medskip
\item{[13]} Fitzsimmons, P.J. and Getoor, R.K.: On the distribution of the
Hilbert transform of the local time of a symmetric L\'evy process. {\it
Ann. Probab.} {\bf 20} (1992), 1484--1497.

\medskip
\item{[14]} Gradshteyn, I.S. and Ryzhik, I.M.: {\it Table of Integrals,
Series, and Products.} Sixth ed. Academic Press, San Diego, CA, 2000.

\medskip
\item{[15]} Grill, K.: On the last zero of a Wiener process. In: {\it
Mathematical Statistics and Probability Theory} (M.L. Puri et al., eds.)
Vol. A, pp. 99--104. D. Reidel, Dordrecht, 1987.

\medskip
\item{[16]} Hu, Y.: The laws of Chung and Hirsch for Cauchy's principal
values related to Brownian local times. {\it Electronic J. Probab.} {\bf
5} (2000), Paper No. 10, 1--16.

\medskip
\item{[17]} Hu, Y. and Shi, Z.: An iterated logarithm law for
Cauchy's principal value of Brownian local times. In: {\it Exponential
Functionals and Principal Values Related to Brownian Motion} (M. Yor,
ed.), pp. 131--154. Biblioteca de la Revista Matem\'atica Iberoamericana,
Madrid, 1997.

\medskip
\item{[18]} Strassen, V.: An invariance principle for the law of the
iterated logarithm. {\it Z. Wahrsch. verw. Gebiete} {\bf 3} (1964),
211--226.

\medskip
\item{[19]} Wen, Jiwei: Some results on lag increments of the principal
value of Brownian local time. {\it Appl. Math. J. Chinese Univ. Ser. B}
{\bf 17} (2002), 199--207.

\medskip
\item{[20]} Yamada, T.: Principal values of Brownian local times and their
related topics. In: {\it It\^o's Stochastic Calculus and
Probability Theory} (N. Ikeda et al., eds.), pp. 413--422. Springer,
Tokyo, 1996.

\medskip
\item{[21]} Yor, M.: {\it Some Aspects of Brownian Motion. Part 1: Some
Special Functionals.} ETH Z\"urich Lectures in Mathematics. Birkh\"auser,
Basel, 1992.

\medskip
\item{[22]} Yor, M., editor: {\it Exponential Functionals and
Principal Values Related to Brownian Motion}. Biblioteca de la Revista
Matem\'atica Iberoamericana, Madrid, 1997.

\bye